\documentclass[11pt,twoside]{article}
\usepackage{graphicx}
\usepackage[sf,bf,compact,topmarks,calcwidth,pagestyles,clearempty]{titlesec}
\usepackage{titletoc}
\usepackage{amsmath}
\usepackage{amsfonts}
\usepackage{indentfirst}
\usepackage{amsthm,amssymb}
\usepackage[nosort]{cite}
\usepackage{color}
\usepackage{xcolor}%
\usepackage{makeidx}
\usepackage[active]{srcltx}

\makeatletter
\renewcommand{\ps@myheadings}{
 \renewcommand{\@oddhead}{
  \parbox{\textwidth}{\hfill\textsl{\rightmark}\hfill
  {\tiny \today} \\
  \rule[.4cm]{\textwidth}{0.5pt}
  }}

\renewcommand{\@evenhead}{
\parbox{\textwidth}{{\tiny \today}\hfill\textsl{\leftmark}\hspace*{\fill}
\\
   \rule[.4cm]{\textwidth}{0.5pt}
  }}
}

\def\R {{\rm I}\hskip -0.85mm{\rm R}}
\def\N {{\rm I}\hskip -0.85mm{\rm N}}

\def\D{\Delta}

\def\G{\Gamma}

\def\m{\mu}

\def\l{\lambda}

\def\O{\Omega}
\def\p{\partial}

\def\n{\nabla}
\def\go{\mathcal}
\def\ds{\displaystyle}
\newtheorem{definition}{Definition}[section]
\newtheorem{theorem}[definition]{Theorem}
\newtheorem{lemma}[definition]{Lemma}

\newtheorem{remark}[definition]{Remark}

\newtheorem{teorema}[definition]{Theorem}

\newtheorem{lema}[definition]{Lemma}

\newtheorem{nota}[definition]{Remark}

\newcommand{\bas}{\begin{eqnarray*}}            %
\newcommand{\eas}{\end{eqnarray*}}              %
\def\qed{{\unskip\nobreak\hfil\penalty50
          \hskip2em\hbox{}\nobreak\hfil\mbox{\rule{1ex}{1ex} \qquad}
   \parfillskip=0pt
   \finalhyphendemerits=0\par\medskip}}
 
 \makeatletter
\begin{document}
\vskip 1cm
\begin{center}
{\bf\Large Long-time behaviour of an angiogenesis model with flux at the tumor boundary}
\end{center}
\begin{center}
{\sc T. Cie\'slak and C. Morales-Rodrigo}
\end{center}
\begin{center}
Institute of Mathematics, Polish Academy of Sciences, 
Warsaw, \'Sniadeckich 8\\
E-mail address: T.Cieslak@impan.pl

Dpto. de Ecuaciones
Diferenciales y An\'alisis Num\'erico\\ Fac. de Matem\'aticas, Univ. de Sevilla,
Calle Tarfia s/n, 41012-Sevilla, Spain\\
E-mail address: cristianm@us.es
\end{center}
\begin{center}
\vskip 0.2cm {\large \it \today} \vskip 0.3cm
\end{center}
\noindent \rule{\textwidth}{.5pt}

\vspace{0.3cm}\noindent{\bf Abstract} \vspace{0.2cm}
\indent

This paper deals with a nonlinear system of partial differential equations modeling a simplified tumor-induced angiogenesis taking into account only the interplay between tumor angiogenic factors and endothelial cells. Considered model assumes a nonlinear flux at the tumor boundary and a nonlinear chemotactic response. It is proved that the choice of some key parameters influences the long-time behaviour of the system. More precisely, we show the convergence of solutions to different semi-trivial stationary states for different range of parameters.

\vspace{0.3cm} \noindent{\it AMS Classification.} 35K45, 35K57, 92C17.
\\
\noindent \rule{\textwidth}{.5pt}

\section{Introduction}
Angiogenesis is a physiological process involving the new vessels sprout from a pre-existing vasculature in response to a chemical stimuli. Angiogenesis is an important ingredient of a processes like development, growth and wound healing. However, angiogenesis is also induced by tumoral cells. In this paper we consider a model of tumor-induced angiogenesis that was proposed in \cite{Mioangioflux}. Actually, in the above mentioned model some factors influencing angiogenesis are neglected to keep the model simple but sufficiently interesting from the analytical point of view. We refer the reader to \cite{mwo} as a source of information about the progress in mathematical modelling and biological knowledge of angiogenesis process. We focus our attention on two key variables: the endothelial cells (ECs), denoted by $u$, and the tumor angiogenic factors (TAF), denoted by $v$. We assume that (ECs) that form the blood vessels wall are induced by the (TAF), factors that are generated by the tumor, to migrate chemotactically towards the tumor. We assume that the (ECs) and the (TAF) fill in a bounded and connected domain $\O \subset \R^d$ with a regular boundary $\p\O$. In particular, neither the existence of extracellular matrix nor the activity of metalloproteinases is considered. But, what was new there, nonlinear flux of TAF on the tumor boundary was taken into account. The reason was that since ECs are supposed to react chemotactically to the TAFs, generating the 
large gradient of TAFs on the boundary would probably make the tumour more dangerous. The aim of \cite{Mioangioflux} was to study the interplay between the density of ECs and TAFs dependently on a parameter $\mu$ measuring the strength of the flux on the tumor boundary and the nonlinearity $V$ measuring nonlinear response of ECs. In \cite{Mioangioflux} the qualitative features of a model were studied in a local sense. We mean by that the local stability of steady states 
which were proven to exist in \cite{Mioangioflux}. We complete the studies taken in \cite{Mioangioflux} by analyzing the global stability of steady states. We shall prove the asymptotic convergence of solutions for different values of $\mu$.
To be more precise, we consider the case
$$
\p\O=\G_1\cup\G_2,
$$
where $\G_1\cap\G_2=\emptyset$ and  $\G_i$ are closed and open sets in the
relative topology of $\p\O$. We suppose that $\G_2$ is the tumor
boundary and $\G_1$ is the blood vessel boundary. Our parabolic problem reads.
\begin{equation} \label{sis00} \left\{
\begin{array}{ll}
u_t-\D u=- \mbox{div}(V (u)\n v)+\l u-u^2 &\mbox{in $\O\times (0,T)$,}\\
v_t-\D v=  -v-cuv &\mbox{in $\O\times (0,T)$,}\\
\displaystyle\frac{\p u}{\p n}= \frac{\p v}{\p n}= 0 &\mbox{on $\G_1\times (0,T)$,}\\
\noalign{\smallskip} \displaystyle\frac{\p u}{\p n}=0,\quad
\frac{\p v}{\p n}=
\m\frac{v}{1+v} &\mbox{on $\G_2\times (0,T)$,}\\
u(x,0)=u_0(x),\quad v(x,0)=v_0(x) &\mbox{in $\O$,}
\end{array}\right.
\end{equation}
where $0<T\leq +\infty$, $\l,\m \in\R$, $c>0$,
\begin{equation}\label{hipotesisV}
V\in \go C^1(\R),\quad V> 0\,\,{\rm in}\,\, (0,\infty)\,\,{\rm
with}\,\, V(0)=0;
\end{equation}
and $u_0$ and $v_0$ are given non-negative and non-trivial 
functions.
In \cite[Theorem 3.1, Theorem 3.8]{Mioangioflux} the existence and uniqueness of global-in-time bounded regular solutions, provided initial data are nonnegative and $V \in L^\infty (0,+\infty)$
is shown. Moreover in \cite[Section 4]{Mioangioflux} the existence of two semi-trivial steady-states $(\lambda,0)$, $\lambda >0$ and $(0,\theta_\mu)$ is shown provided $\mu > \mu_1$ (see also \cite{UmezuNoDEA}), where $\mu_1$ is the principal eigenvalue of the boundary eigenvalue problem
$$
\left \{ \begin{array}{ll}
-\Delta v + v =0 & \mbox{in $\Omega$,}\\
\displaystyle\frac{\p v}{\p n} =0 & \mbox{on $\Gamma_1$,}\\[2mm]
\displaystyle\frac{\p v}{\p n} =\mu v & \mbox{on $\Gamma_2$.} 
\end{array}
\right.
$$
Furthermore, results concerning the linearized stability around the semi-trivial solutions to (\ref{sis00}) are proven in \cite{Mioangioflux}. \\

First models of tumor induced angiogenesis that we are aware of are considered in \cite{AnCha} (see also \cite{LevineJMB} for a more elaborated model). A reduced model proposed in \cite{LevineJMB} is studied in \cite{FontelosSIAM}. The local stability of the homogeneous steady-states in one dimensional domains is shown there. In all the mentioned papers the boundary conditions are either zero Neumann or no-flux. In \cite{SuarezJDE2008} the stationary problem of (\ref{sis00}) with linear flux for $v$ is studied. Finally let us mention \cite{Kettemann} where the authors study the local solvability of a system of partial differential equations with a nonlinear boundary condition and a chemotaxis term.\\
 
The aim of this paper is to analyze the global stability for positive initial data. In particular we show global stability for some range of parameters $(\lambda,\mu)$ for which even the local stability is not known.\\

It should be pointed out that the results of this paper could be extended even to more general forms of $V$ as soon as
\begin{equation}\label{unbreak}
\|u(t)\|_\infty <C, \;\;\mbox{for}\;\; t\geq 0.
\end{equation}
Observe that, if the above inequality holds, then the parabolic regularity asserts $$\|v(t)\|_\infty <C$$ for any $t>0$ and by \cite[Theorem 3.1]{Mioangioflux} the solution is global and regular. In particular, when $V$ is bounded in the $L^\infty$ norm (see \cite{Mioangioflux}) then (\ref{unbreak}) is satisfied.
 \section{Preliminaries}
For the reader's convenience we collect here some results of interpolation theory and its applications to parabolic problems that will be used throughout the paper.
\begin{enumerate}
\item Let $E_0$, $E_1$ two normed spaces, we can define the real interpolation functor, denoted by 
$$
(E_0,E_1)_{\theta,p}, \;\; 0<\theta<1, \;1\leq p \leq +\infty,
$$
(see for instance \cite[Def. 22.1]{Tartar}). During the paper we will use the following property of the real interpolation functor (see \cite[Lemma 25.2]{Tartar}):\\
If $(E_0,E_1)_{\theta,p}$ is a Banach space then
\[
\exists C >0\;\; \mbox{such that} \;\; \|x\|_{(E_0,E_1)_{\theta,p}} \leq C \| x\|_{E_0}^{1-\theta} \|x\|_{E_1}^\theta \;\; \forall x \in E_0 \cap E_1.
\]
In the context of fractional Sobolev spaces this inequality reads, see \cite[Theorem 7.2]{amann1993}
\begin{equation}\label{interpolation}
\| u(t)\|_{W^{m,p}} \leq C \|u(t)\|_{W^{k,p}}^\theta
\| u(t)\|_{p}^{1-\theta}
\end{equation}
for $m<k\theta$, $\theta\in (0,1)$.
\item Let us consider a parabolic problem with a non-homogeneous boundary condition 
\begin{equation}\label{linear-p}
\left \{ \begin{array}{ll}
z_t + \mathcal{A} z = f(t) &\mbox{ in $\Omega\times (0,T)$,}\\
\mathcal{B}z = g(t) &\mbox{ on $\partial\Omega \times (0,T)$,}\\
z(x,0)=z_0 (x), &\mbox{ in $\Omega$}.
\end{array}
\right.
\end{equation}
where 
$$
\mathcal{B} z := \frac{\partial z}{\partial n}
$$
and 
$$
\mathcal{A}z: = -\Delta z + z.
$$
We define the space of functions  
$$
W^{s,p}_{\go{B}} := \left \{ \begin{array}{ll}
\{ z \in W^{s,p} (\Omega): \go{B} z =0 \} & \mbox{ if $1+1/p<s\leq 2$},\\
W^{s,p} (\Omega) & \mbox{ if $-1+1/p <s<1+1/p$},\\
(W^{-s,p'}(\Omega))' & \mbox{ if $-2+1/p <s\leq -1+1/p$}.
\end{array}
\right.
$$
It is known that $(\go{A},\go{B})$, as being in separated divergence form (see  \cite[pg. 21]{amann1993}), is normally elliptic. We denote by $A_{\alpha -1}$ the $W^{2\alpha -2,p}_{\go{B}}$-realization of $(\go{A},\go{B})$ (see \cite[pg. 39]{amann1993} for the precise definition). Since $(\go{A},\go{B})$ is normally elliptic then $A_{\alpha-1}$ generates an analytic semigroup \cite[Theorem 8.5]{amann1993}. Moreover, if 
$$
(f,g) \in \go C ((0,T);W^{2 \alpha -2,p}_{\go B} (\Omega) \times W^{2\alpha -1 -1/p,p}_{\go B} (\partial \Omega) )
$$
for some $T>0$ and $2\alpha \in (1/p, 1+1/p)$ then for any $t<T$ we rewrite (\ref{linear-p}) by the generalized variation of constants formula
$$
z(t)= e^{-t A_{\alpha-1} } z_0 + \int_0^t e^{-(t-\tau)A_{\alpha-1}} (f(\tau)+A_{\alpha -1} {\go B}^c_{\alpha} g (\tau)) d\tau, 
$$
where $\go{B}^c_{\alpha}$ is the continuous extension of $(B |_{Ker(\go A)})^{-1}$ to $W^{2\alpha -1-1/p,p} (\partial \Omega)$. Since $[0,+\infty) \subset  \rho (-A_{\alpha-1})$ ($\rho$ is the resolvent set) then by \cite[Remark 8.6 c)]{amann1993} there exists a constant $C \geq 1$ such that
\begin{equation}\label{P-Aalpha}
\|z\|_{W^{2\alpha,p}_{\go B} } \leq C \|A_{\alpha-1} z\|_{W^{2\alpha-2,p}_{\go B}}.
\end{equation}
\item Let $a,b,c \in L^\infty (\Omega)$, the eigenvalue problem
$$
\left \{ \begin{array}{ll}
-\Delta z + a(x) z = \lambda z & \mbox{ in $\Omega$},\\[2mm]
\displaystyle\frac{\partial z}{\partial n} + b(x) z =0 & \mbox{ on $\Gamma_1$},\\[2mm]
\displaystyle\frac{\partial z}{\partial n} + c(x)z =0 & \mbox{ on $\Gamma_2$}.
\end{array}
\right.
$$
has a unique principal eigenvalue (i.e. an eigenvalue whose associated eigenfunction can be chosen positive in $\Omega$) and it will be denoted by
$$
\lambda_1 (-\Delta + a;{\go N}+ b; {\go N}+c ).
$$
\end{enumerate}
\section{Convergence to the semi-trivial solution $(\l,0)$}
In the present section we deal with the convergence to the semi-trivial steady-state $(\l,0)$. Throughout this section we assume (\ref{unbreak}). A sufficient condition guaranteeing (\ref{unbreak}) is the boundedness of $V$ (see \cite{Mioangioflux}). We will use the generalized variation of constants formula to estimate $v$, which is stated in the next lemma. 
\begin{lema} \label{Est-v-Aalhpa}
Let $\gamma \in (1,+\infty)$ and $\beta \in (1, 2\alpha)$. Then, for every $\tau \in (0,t)$ there exists a constant $\delta \in (0,Re\; \sigma (A_{\alpha-1}))$ ($\sigma$ denotes the spectrum)  and $\theta=\theta (\beta) \in (0,1)$ such that
$$
\| e^{-(t-\tau)A_{\alpha-1}} z \|_{W^{\beta,\gamma}} \leq C (t-\tau)^{-\theta}e^{-\delta (t -\tau)} \|z\|_{W^{2\alpha-2,\gamma}_{\go B}}
$$
for every $z \in W^{2\alpha,\gamma}_{\go B}$.
\end{lema}
\noindent {\bf Proof.} By the choice of $\beta$ we have $W^{\beta, \gamma}_{\go B} = W^{\beta,\gamma} (\Omega)$. As a consequence if we apply  \cite[Theorem 7.2]{amann1993} we get
$$
\| e^{-(t-\tau)A_{\alpha-1}} z \|_{W^{\beta,\gamma}} \leq C \| e^{-(t-\tau)A_{\alpha-1}} z \|_{W^{2\alpha,\gamma}_{\go B}} \|^\theta \|e^{-(t-\tau)A_{\alpha-1}} z \|_{W^{2\alpha -2,\gamma}_{\go B}}^{1-\theta} 
$$
for some $\theta \in (0,1)$. Next we apply (\ref{P-Aalpha}) to the first norm on the right hand side and \cite[Theorem 1.3.4]{henry1981}
to deduce 
$$
\|  e^{-(t-\tau)A_{\alpha-1}} z \|_{W^{\beta,\gamma}} \leq C (t-\tau)^{-\theta} e^{-\delta (t-\tau)\theta }e^{-\delta (t-\tau)(1-\theta)}\| z\|_{W^{2\alpha-2,\gamma}_{\go B}},
$$
where $\delta \in (0,Re\; \sigma (A_{\alpha-1}))$. \qed
\begin{lema}\label{convergv}
Let $\gamma \in (1,+\infty)$, $\beta \in (1,1+1/\gamma)$, $\mu \in [0,\mu_1)$ and $0<\delta <\rho<\alpha(\mu)$ where $\alpha (\mu)$ is defined as
$$
\alpha(\mu):= \lambda_1 (-\Delta+1;\go N,\go N-\mu).
$$
Then, there exists $C>0$ such that, for $t>0$, the $v$-solution to
(\ref{sis00}) satisfies
$$
\begin{array}{ll}
v(x,t)&\leq Ce^{-\rho t} \;\;\;\; \forall (x,t) \in \overline \Omega \times (0,+\infty),\\
\|v(t)\|_{W^{\beta,\gamma}} &\leq C (1+t^{-\theta})e^{-\delta t}\|v_0\|_{\gamma},
\end{array}
$$
where $\theta = \theta (\beta ) \in (0,1)$.
\end{lema}
\noindent {\bf Proof.} A solution to the problem
\begin{equation}\label{eqw} \left \{
\begin{array}{ll}
w_t -\Delta w +w = 0 & \mbox{ in $\Omega \times (0,T_{max})$},\\
\ds\frac{\partial w}{\partial n} = 0 & \mbox{ on $\Gamma_1 \times (0,T_{max})$},\\[2mm]
\ds\frac{\partial w}{\partial n} = \mu w & \mbox{ on $\Gamma_2 \times (0,T_{max})$},\\
w(x,0)=v_0 (x) & \mbox{ in $\Omega$,}
\end{array}
\right.
\end{equation}
is a supersolution to the $v$-equation of (\ref{sis00}), therefore
$v(x,t) \leq w(x,t)$. Since, for sufficiently large $M$, $\overline w=M e^{-\rho t} \varphi_1$, with $\varphi_1$ a positive eigenfunction associated to $\alpha (\mu)$, is a supersolution to (\ref{eqw}), the pointwise estimate in the claim of the lemma follows. For the second one we pick
$$
f(t):= -cu(t)v(t),
$$
$$
g(t):= \left \{ \begin{array}{ll}
0 & \mbox{ on $\Gamma_1$},\\
\displaystyle \mu \frac{v(t)}{1+v(t)} & \mbox{ on $\Gamma_2$}\;.
\end{array}
\right.
$$
Taking the $W^{\beta,\gamma}$-norm in a generalized variation of constants formula for $v$ and using Lemma \ref{Est-v-Aalhpa} we obtain
$$
\begin{array}{ll}
\| v(t) \|_{W^{\beta,\gamma}} \hspace{-10pt} & \leq \displaystyle \|e^{-t A_{\alpha-1}} v_0\|_{W^{\beta,\gamma}} +\int_0^t \|e^{-(t-\tau)A_{\alpha-1}}(f(\tau)+A_{\alpha-1} {\go B}^c_\alpha g(\tau))\|_{W^{\beta,\gamma}}\\
&\leq \displaystyle C \left ( e^{-\delta t } t^{-\theta} \|v_0\|_{W^{2\alpha-2,\gamma}_{\go B}} +\int_0^t (t-\tau)^{-\theta}e^{-\delta (t-\tau)} \|f(\tau)+A_{\alpha-1} {\go B}^c_{\alpha} g (\tau)\|_{W^{2\alpha-2,\gamma}_{\go B}} d\tau \right ).
\end{array}
$$
Next, we estimate the last term in the above inequality using the fact that
$$A_{\alpha-1} {\go B}^c_{\alpha}  \in {\go L}(W^{2\alpha-1-1/\gamma,\gamma} (\partial \Omega), W^{2\alpha-2,\gamma}_{\go B} (\Omega))$$ 
and the continuous embeddings 
$$
L^\gamma (\Omega) \hookrightarrow W^{2\alpha-2,\gamma}_{\go B},\qquad L^\gamma (\partial \Omega) \hookrightarrow W^{2\alpha-1-1/\gamma,\gamma} (\partial \Omega).
$$
Therefore, we get 
\begin{equation}\label{convw1p}
\|v(t)\|_{W^{\beta,\gamma}} \leq e^{-\delta t} t^{-\theta} \|v_0\|_{\gamma} + Ce^{-\delta t} \int_0^t e^{\delta \tau} (t-\tau)^{-\theta} \left(\|f(\tau)\|_{L^\gamma (\Omega)} + \|g(\tau)\|_{L^\gamma (\partial \Omega)} \right) d\tau.
\end{equation}
Observe that by (\ref{unbreak}) and the first part of the Lemma we have
$$
\| f(\tau) \|_{L^\gamma (\Omega)} \leq C \|v\|_{L^\infty(\Omega)} \leq Ce^{-\rho \tau} ,
$$
$$
\|g(\tau)\|_{L^\gamma (\partial \Omega)} \leq \|v\|_{L^\infty (\partial \Omega)} \leq Ce^{-\rho \tau} .
$$
In view of the above bounds, (\ref{convw1p}) yields
$$
\|v(t) \|_{W^{\beta,\gamma}} \leq C e^{-\delta t} t^{-\theta} \|v_0\|_{\gamma} + C e^{-\delta t} \|v_0\|_{\gamma} \int_0^t e^{(\delta-\rho)\tau} (t-\tau)^{-\theta} d\tau.
$$
Next, by the choice of $\delta$ and $\rho$, $\displaystyle{\int_0^{\infty} e^{(\delta-\rho) \tau} (t -\tau)^{-\theta} d\tau = C < +\infty}$ and the Lemma follows.
\qed

Our purpose is to show that $u$ converges to steady states.
To this end we treat separately the cases $\lambda=0$, $\lambda>0$.\\
\subsection{Case $\lambda=0$.}
\begin{lema}\label{convLema}
Let $\tau>0$ and $y \in C^1(\tau,+\infty) \cap L^1
(\tau,+\infty)$, $y'\in L^1(\tau,+\infty)$. 
Then $\displaystyle{\lim_{t\rightarrow +\infty} |y(t)|=0}$.
\end{lema}
\noindent {\bf Proof.} By the assumptions of the lemma we observe that for any $k>0$
\begin{equation}\label{ten}
\lim_{t\rightarrow +\infty} \int_t^{t+k} \left(|y(s)|+ |y'(s)| \right)ds =0.
\end{equation}
Let us assume that $\displaystyle{\lim_{t\rightarrow +\infty} |y(t)|\not = 0}$, then there exists a sequence $\{ t_n\}_{n\in\N}$, $t_n \rightarrow +\infty$, such that
$$
|y(t_n)| >C>0, 
$$
for all $n \geq n_0$.
We pick $\theta \in (0,k]$, then for any $\varepsilon>0$
$$
\big ||y(t_n+\theta)|-|y(t_n)|\big |\leq |y(t_n+\theta) - y(t_n)|\leq
\int_{t_n}^{t_n + \theta} |y'(s)|ds \leq \int_{t_n}^{t_n +k}
|y'(s)|ds<\varepsilon
$$
by (\ref{ten}).
Therefore $|y(s)| > C/2$ for all $s\in [t_n,t_n+k]$, $n\geq n_0$.
The last assertion contradicts the fact that
$$
\lim_{n\rightarrow +\infty} \int_{t_n}^{t_n+k} |y(s)|ds =0.
$$
\qed
In the following lemmata $(u,v)$ is a solution to (\ref{sis00}).
\begin{lemma} \label{lforu} Let $\lambda =0$ and $t>\tau>0$, then it holds
\begin{equation}\label{ss1}
\displaystyle{ \mu \int_\tau^t \int_{\Gamma_2} \frac{V(u)v}{1+v} +
\int_\tau^t \int_\Omega u^2 = \int_\Omega u(\tau) -\int_\Omega
u(t).}
\end{equation}
\end{lemma}
\noindent {\bf Proof.} Integrating the $u$-equation of
(\ref{sis00}) yields
$$
\begin{array}{ll}
\displaystyle{\int_\Omega u_t} &= \displaystyle{ \int_{\partial \Omega} \left ( \frac{\partial u}{\partial n} - V(u)\frac{\partial v}{\partial n}\right )- \int_\Omega u^2 } \\[2mm]
&= \displaystyle{-\mu \int_{\Gamma_2} \frac{V(u) v}{1+v}} -
\int_\Omega u^2.
\end{array}
$$
So, integrating the last expression in time between
$\tau$ and $t$ we get the result. \qed
\begin{remark} \label{i-remark}
By Lemma \ref{lforu} we see that for any $t>\tau$
$$
\int_\tau^t \int_\Omega u^2 \leq \|u(\tau)\|_1. 
$$
\end{remark}
\begin{teorema} \label{convul0}
Assume that $0\leq \mu <\mu_1$ and $\lambda=0$, then
$$
\lim_{t\rightarrow +\infty} \|u(t)\|_{W^{m,p}}
=0,
$$
for any $m<1$ and $p\geq 2$.
\end{teorema}
\noindent {\bf Proof.} On multiplying the $u$-equation of (\ref{sis00})
by $u$ and integrating in space we obtain 
\begin{equation} 
\begin{array}{ll}
\displaystyle{ \frac{d}{2dt} \int_\Omega u^2 } &=\displaystyle{ \int_\Omega \left ( -|\nabla u|^2 + V(u)\nabla v \cdot \nabla u -u^3\right )- \mu\int_{\Gamma_2} \frac{V(u)uv}{1+v}} \\
&\leq \displaystyle{(\epsilon -1) \int_{\Omega}|\nabla u|^2 +
C(\epsilon)\int_{\Omega} |\nabla v|^2 - \mu \int_{\Gamma_2}
\frac{V(u)uv}{1+v}-\int_{\Omega} u^3.}
\end{array}
\end{equation}
Therefore, we infer
$$
\displaystyle{ \frac{d}{2dt} \int_\Omega u^2 +
(1-\epsilon)\int_{\Omega} |\nabla u|^2 \leq C(\epsilon)
\|v\|_{W^{1,2}}^2,}
$$
and after integrating in time, thanks to Lemma \ref{convergv} we
arrive at
$$
\displaystyle{\int_{\Omega} u(t)^2 - \int_\Omega u(\tau)^2 +
(1-\epsilon) \int_\tau^t \int_\Omega |\nabla u|^2 \leq C(\epsilon)
\int_\tau^t (1+s^{-\theta})^2 e^{-2\delta s} \|v_0\|_{p}^2}.
$$
In particular we deduce that for $t>\tau$
$$
\displaystyle{\int_\tau^t \int_\Omega |\nabla u|^2 \leq C .}
$$
By \cite[Lemma 3.8]{Mioangioflux} we find a bound $\|u(t)\|_{C(\overline \Omega)} \leq C$, therefore, 
$$
\displaystyle{ \left |\frac{d}{2dt} \int_{\Omega} u^2 \right |
\leq C\int_{\Omega} |\nabla u|^2 + C(\epsilon) \| v\|_{W^{1,2}}^2
+ C \mu \int_{\Gamma_2} \frac{ V(u)v}{1+v}+ C\int_\Omega u^2}.
$$
Thanks to (\ref{ss1}), for $t >\tau$
\begin{equation}\label{imp2}
\displaystyle{ \int_\tau^t\left |\frac{d}{2dt} \int_{\Omega} u^2
\right | \leq C.}
\end{equation}
Finally, Remark \ref{i-remark} and (\ref{imp2}) together with
Lemma \ref{convLema} entail
\[
\lim_{t\rightarrow +\infty} \|u(t)\|_{2} =0.
\]
Also thanks to $\|u(t)\|_{C(\overline \Omega)} \leq C$ for all
$t>0$ we obtain
\[
\lim_{t\rightarrow +\infty} \|u(t)\|_{p} =0
\]
for any $p>2$. 
Next we recall that by \cite[Lemma 3.7]{Mioangioflux} for any $2\beta \in (k,1)$ we find a bound on the $X_\beta$ norm of $u$, where 
$X_\beta$ is a usual fractional space connected to a semigroup approach to parabolic equations, see \cite{henry1981}. Next, due to the fact that $2\beta \in (k,1)$, we infer from the embedding $X_\beta \hookrightarrow W^{k,p}$ (see for instance \cite[Theorem
1.6.1]{henry1981}) that for all $k<1$ and $p\geq 2$
\[
\| u(t)\|_{W^{k,p}} \leq C.
\]
Next, (\ref{interpolation}) entails
$$
\| u(t)\|_{W^{m,p}} \leq C \|u(t)\|_{W^{k,p}}^\theta
\| u(t)\|_{p}^{1-\theta}.
$$
Therefore, it holds
\begin{equation}
\lim_{t \rightarrow +\infty} \| u(t)\|_{W^{m,p}} \leq C
\lim_{t\rightarrow +\infty} \|u(t)\|_{p}^{1-\theta} =0.
\end{equation}
\qed
\begin{nota}
Let us point out that if we pick $m$ such that $m-d/p >0$ then $W^{m,p} (\Omega)$ is embedded in $\go C(\overline \Omega)$.
\end{nota}
\subsection{Case $\lambda>0$.}
Assume that there exists $\delta_0$ and $t_0$ such that
\begin{equation}\label{H} 
u(t) >\delta_0>0
\end{equation}
for $t>t_0>0$.
Next, we examine the long time behavior for $u$ under the hypothesis
$(\ref{H})$. In the sequel we shall give sufficient conditions on $V(u)$
implying $(\ref{H})$.
\begin{teorema}
Let $0\leq \mu < \mu_1$ and assume the the hypothesis $(\ref{H}))$ is
satisfied, then there exists $\theta >0$ such that
\begin{equation} \label{yang}
\|u(t) -\lambda\|_{W^{m,p}} \leq C e^{-\theta
t},
\end{equation}
for all $t\geq t_0$ and any $m<1$, $p\geq 2$.
\end{teorema}
\noindent {\bf Proof.} On multiplying the $u$-equation by $u-\lambda$ we
have
\begin{equation}
\begin{array}{ll}
\displaystyle{\frac{d}{2dt}\int_\Omega (u-\lambda)^2} = & \displaystyle{-\int_\Omega |\nabla u|^2 + \int_\Omega V(u)\nabla v \cdot \nabla u-\mu \int_{\Gamma_2}\frac{vV(u)}{1+v}(u-\lambda) - \int_\Omega u(u-\lambda)^2}\\
&\hskip -6mm \leq \hskip 2mm \displaystyle{-\frac{1}{2} \int_\Omega |\nabla u|^2 + \frac{\|V\|_{\infty}^2}{2}\int_\Omega |\nabla v|^2+}\\
&+\displaystyle{\mu \|V(u)(u-\lambda)\|_{2,\Gamma_2}
\left(\int_{\Gamma_2} \frac{v^2}{(1+v)^2}\right)^{1/2} -
\int_\Omega u(u-\lambda)^2}.
\end{array}
\end{equation}
Having in mind that $(1+v)^2 \geq 1$, the hypothesis $(\ref{H})$ and the
Sobolev trace embedding $$W^{1,2}(\Omega) \hookrightarrow
L^2(\partial\Omega)$$ we get
\begin{equation} \label{endc}
\displaystyle{\frac{d}{dt}\int_\Omega (u-\lambda)^2 + 2\delta_0
\int_\Omega (u-\lambda)^2 \leq C \|v\|_{W^{1,2}}^2 + \mu C
\|v\|_{W^{1,2}}}.
\end{equation}
By Lemma \ref{convergv} we can deduce
$$
\|u(t)-\lambda\|_{2}^2 \leq Ce^{-\theta_1 t}
$$
for $0<\theta_1 <\min\{2\delta_0,\beta\}$. At this point we can
argue exactly as in the end of the proof of Theorem~\ref{convul0}. Namely, by the bound on $u$ in $L^\infty$
we infer the bound on the $L^p$ norm of $u$, $p>2$. Next, we use the estimate of $u$ in $W^{k,p}$, $k<1,p\geq 2$,
coming from \cite[Lemma 3.7]{Mioangioflux}, in order to conclude (\ref{yang}).
\qed

In the rest of this section we give sufficient conditions on $V$ implying
$(\ref{H})$. Actually, only the behavior of $V$
around zero matters. Roughly speaking we require a superlinear growth of $V$
in the neighbourhood of zero. From now on we assume that there exist
$C,\delta>0$, $k_0>1+d/2$, $j>d/2$ such that
\begin{equation}\label{H1}
0<V(s)<Cs^{k_0}\,,\hspace{20pt} |V'(s)|\leq Cs^j\,
\end{equation}
for all $s\in (0,\delta)$.
\begin{nota}
The condition $(\ref{H1})$ is satisfied, for example, for functions
$$\displaystyle{V(u) = \frac{u^\alpha}{1+u^\alpha}}$$ with $\alpha
>1+d/2$.
\end{nota}
Next we introduce some notation that will be of importance in the proof of
(\ref{H}). Moreover we formulate a lemma which we need in the main part of the proof of 
(\ref{H}).
Let $f(\delta),g(\delta)$ be defined in a following way:
$$
\displaystyle{f(\delta):=\sup_{s\in (0,\delta)} V^2(s)},
$$  
$$
\displaystyle{ g(\delta):=\sup_{s\in(0,\delta)} (2(s-\delta)_-^2
V'(s)^2 + 2V^2(s)).}
$$
\begin{lema}\label{Lemato}
Assume that (\ref{H1}) holds. Moreover, for some $D,\mu>0$, $\eta>1$, $\tilde{\epsilon}$ and $C(\tilde{\epsilon})$ are given by 
\[
\tilde{\epsilon}=\frac{\delta^{2\eta}}{2\mu D},\;\;C(\tilde{\epsilon})=\frac{\mu D}{2\delta^{2\eta}}.
\] 
Then, if $\delta>0$ is small enough, the following conditons are satisfied simultaneously
\begin{equation}\label{cond2}
C(\widetilde \epsilon) \frac{V^2(s)}{s}\delta \leq \lambda -\delta
\end{equation}
for $s\in(0,\delta)$,
\begin{equation} \label{cond1}
C(\widetilde \epsilon)g(\delta) <1/2
\end{equation}
and
\begin{equation} \label{cond3}
f(\delta)D \leq \frac{\delta^{2\eta}}{2}.
\end{equation}
\end{lema}
\noindent {\bf Proof.} Thanks to $(\ref{H1})$, we have
$$
f(\delta)D=\sup_{s \in (0,\delta)} V^2(s)D
\leq C\delta^{2k_0}D.
$$
Hence, for $\eta <k_0$ and $\delta$ sufficiently small (\ref{cond3}) is satisfied. Next, owing to $(\ref{H1})$, we
observe that
$$
\displaystyle{C(\widetilde \epsilon) \frac{V^2(s)}{s}\delta } \leq C(\widetilde \epsilon) \delta^{2k_0}.
$$
Thus, (\ref{cond2}) can be assured for $\eta <k_0$ and
$\delta$ small enough. Moreover, it is straightforward to see that
(\ref{cond1}) is also satisfied for
$1<\eta<\min\{k_0,1+j\}$.\qed
\begin{lema}\label{Lematon}
Assume that $0\leq \mu <\mu_1$ and that $(\ref{H1})$ is satisfied then
$(\ref{H})$ holds.
\end{lema}
\noindent {\bf Proof.} Let $\delta >0$ be a fixed constant defined in (\ref{H1}). Given a function $f$, we define the negative part of $f$
as a nonpositive function as follows 
\[
f_- :=\min \{ f,0 \}.
\]
Our purpose is to show that $\| (u -\delta)_- (t)\|_{\infty} \leq \delta/2$ for every $t>t_0$ which implies (\ref{H}). In order to obtain the previous estimate we
multiply the $u$-equation by $(u-\delta)_-$ and we integrate in
space to obtain
$$
\begin{array}{ll}
\displaystyle{\frac{d}{2dt} \int_\Omega (u-\delta)_-^2} =& - \displaystyle{\int_\Omega (\nabla u - V(u) \nabla v)\cdot \nabla (u-\delta)_-}\\[2mm]
&+ \displaystyle{\int_{\partial \Omega} \left ( \frac{\partial u}{\partial n}- V(u)\frac{\partial v}{\partial n}\right)(u-\delta)_- +} \displaystyle{\int_\Omega u(\lambda -u)(u-\delta)_-}\\[2mm]
&\hskip -6mm = \hskip 2mm \displaystyle{-\int_\Omega |\nabla (u-\delta)_-|^2 + \int_\Omega V(u)\nabla v \cdot \nabla (u-\delta)_-}\\[2mm]
&-\displaystyle{\int_{\Gamma_2} V(u) \mu \frac{v}{1+v} (u-\delta)_- + \int_\Omega u(\lambda-u)(u-\delta)_-}\\
& \hskip -6mm = \hskip 2mm \displaystyle{-\int_\Omega |\nabla (u-\delta)_-|^2 + \int_{\Omega_\delta} V(u) \nabla v \cdot \nabla (u-\delta)_-}\\[2mm]
&- \mu \displaystyle{ \int_{\Gamma_\delta} \frac{v}{1+v}
V(u)(u-\delta)_-+ \int_\Omega u(\lambda -u)(u-\delta)_-},
\end{array}
$$
where
$$
\Omega_\delta := \{x \in \Omega :\; u(x)<\delta \}\,,\;\;\;
\Gamma_\delta:=\{ x \in \Gamma_2:\; u(x)<\delta\}.
$$
Consequently,
$$
\begin{array}{ll}
\displaystyle{\frac{d}{2dt}\int_\Omega (u-\delta)_-^2} \leq & \displaystyle{(\epsilon -1)\int_\Omega |\nabla (u-\delta)_-|^2 +C(\epsilon) \int_{\Omega_\delta} V^2(u) |\nabla v|^2}\\
&- \displaystyle{\mu \int_{\Gamma_\delta} \frac{v}{1+v} V(u)(u-\delta)_- + \int_\Omega u(\lambda -u)(u-\delta)_-}\\
& \hskip -6mm \leq \hskip 2mm \displaystyle{(\epsilon -1)\int_\Omega |\nabla (u-\delta)_-|^2 +C(\epsilon) \sup_{s\in (0,\delta)} V^2(s) \int_\Omega |\nabla v|^2}\\
&- \displaystyle{\mu \int_{\Gamma_\delta} \frac{v}{1+v}
V(u)(u-\delta)_- + \int_\Omega u(\lambda -u)(u-\delta)_- .}
\end{array}
$$
Previous inequality can be rewritten in terms of $f(\delta)$ defined before Lemma \ref{Lemato} as
$$
\begin{array}{ll}
\displaystyle{\frac{d}{2dt}\int_\Omega (u-\delta)_-^2} \leq & \displaystyle{(\epsilon -1)\int_\Omega |\nabla (u-\delta)_-|^2 +C(\epsilon)f(\delta) \int_\Omega |\nabla v|^2 +\mu \widetilde \epsilon \int_{\Gamma_2} \frac{v^2}{(1+v)^2}}\\
&+ \displaystyle{ \mu C(\widetilde \epsilon) \int_{\Gamma_2}
V(u)^2(u-\delta)_-^2 + \int_\Omega u(\lambda-u)(u-\delta)_-}.
\end{array}
$$
Thanks to the Sobolev trace embedding $W^{1,2} (\Omega )
\hookrightarrow L^2(\partial \Omega)$ and having in mind that
$(v+1)^2 \geq 1$, we arrive at
$$
\begin{array}{ll}
\displaystyle \int_{\Gamma_2} V(u)^2 (u-\delta)_-^2 \leq  & \displaystyle C \left(
\int_\Omega V^2(u)(u-\delta)_-^2 + \right. \\ 
& \hspace{15pt} \left .  \displaystyle +\int_\Omega \left(2(u-\delta)_-^2
V'(u)^2+2 V^2(u)\right) |\nabla (u-\delta)_-|^2\right),
\end{array}
$$
$$
\displaystyle{\mu \widetilde \epsilon \int_{\Gamma_2}
\frac{v^2}{(1+v)^2} \leq C \mu \widetilde \epsilon \|v\|_{W^{1,2}}^2.}
$$
Therefore, we obtain
\[
\frac{d}{2dt}\int_\Omega (u-\delta)_-^2 \leq  (\epsilon -1)\int_\Omega |\nabla (u-\delta)_-|^2 +C(\epsilon)f(\delta) \int_\Omega |\nabla v|^2 + C\widetilde \epsilon \|v\|_{W^{1,2}}^2
\]
\[
+C(\widetilde\epsilon)\left(\int_\Omega V^2(u) (u-\delta)_-^2 +\int_\Omega \left(2(u-\delta)_-^2 V'(u)^2+2 V^2 (u)\right)|\nabla (u-\delta)_-|^2\right)
\]
\begin{equation}\label{wzor}
+\int_\Omega u(\lambda -u)(u-\delta)_-.
\end{equation}
In view of the nonnegativity of $u$ we have
\begin{equation}\label{numer}
-\delta<(u-\delta)_{-}.
\end{equation}
Owing to (\ref{numer}), from (\ref{wzor}) we see that ($g(\delta)$ was defined before Lemma \ref{Lemato})
\[
\frac{d}{2dt}\int_\Omega (u-\delta)_-^2\leq (\epsilon +C(\widetilde \epsilon) g(\delta)-1 )\int_\Omega |\nabla (u-\delta)_-|^2
\]
\begin{equation}\label{wzor1}
+(C(\epsilon)f(\delta) + \mu \widetilde \epsilon)\|v\|_{W^{1,2}}^2
+\int_\Omega u(u-\delta)_-\left ( \lambda - u -
C(\widetilde \epsilon) \frac{V^2(u)}{u} \delta \right).
\end{equation}
Due to the nonpositivity of $(u-\delta)_{-}$ and (\ref{cond1}) we have
\begin{equation}\label{wzor2}
\int_\Omega u(u-\delta)_-\left ( \lambda - u -C(\widetilde \epsilon) \frac{V^2(u)}{u} \delta \right)<0.
\end{equation}
By the Hopf lemma and zero Neumann data on the boundary for $u$ we see that there exists $\delta_1$ such 
that $u(t_0)>\delta_1$. Hence choosing $\delta<\delta_1$ and using (\ref{cond2}), (\ref{wzor2}) and Lemma \ref{convergv} we infer from (\ref{wzor1})  
$$
\|(u-\delta)_- (t)\|_{2}^2 \leq (2C(\epsilon)f(\delta) + 2\mu
\widetilde \epsilon)C(\beta),
$$
for $t >t_0 >0$.
We shall show that
\begin{equation}\label{L2}
\|(u-\delta)_-(t)\|_{2}^2 \leq \delta^{2\eta},
\end{equation}
for some $\eta>1$. To this
end notice that choosing $\epsilon=C(\epsilon)=1/2$, we are in a position to apply Lemma \ref{Lemato} with 
$D=C(\beta)$. As a consequence, for $\tilde{\epsilon}$ as it is chosen in Lemma \ref{Lemato},  
(\ref{cond1}),(\ref{cond2}), (\ref{cond3}) and 
\[
2\mu\tilde{\epsilon} C(\beta)\leq\frac{\delta^{2\eta}}{2}
\] 
are satisfied simultaneously. Hence (\ref{L2}) is shown.

Next we use interpolation between $L^p$ spaces, (\ref{L2}) and (\ref{numer}) to obtain
$$
\begin{array}{ll}
\|(u-\delta)_-\|_{2/\theta_1} &\leq \|(u-\delta)_-\|_{2}^{\theta_1} \|(u-\delta)_-\|_{\infty}^{1-\theta_1}\\
&\leq \delta^{\alpha \theta_1}
\delta^{1-\theta_1}=\delta^{1+(\alpha-1)\theta_1}.
\end{array}
$$
Applying (\ref{interpolation}), we infer
$$
\|(u-\delta)_-\|_{W^{\theta,2/\theta_1}} \leq C
\|(u-\delta)_-\|_{W^{1,2/\theta_1}}^\theta
\|(u-\delta)_-\|_{2/\theta_1}^{1-\theta}\leq C_1\|(u-\delta)_-\|_{2/\theta_1}^{1-\theta},
$$
the last inequality being a consequence of the uniform bound
of $L^\infty$ norm, see \cite[Theorem Lemma 3.8]{Mioangioflux}, and \cite[Theorem
15.5]{amann1993}.
Picking up $\theta_1$ such that
\begin{equation} \label{cond5}
\theta -\frac{d\theta_1}{2}>0
\end{equation}
we make sure that $W^{\theta,2/\theta_1} (\Omega) \hookrightarrow L^\infty (\Omega)$.
Consequently,
\[
\|(u-\delta)_-\|_\infty\leq C_2\delta^{(1-\theta)(1+(\alpha-1)\theta_1)}.
\]
Next, we notice that choosing $\alpha>1+\frac{d}{2}$ we make sure that
\[
1< \left( 1-\frac{d\theta_1}{2} \right) (1+(\alpha-1)\theta_1).
\]
Hence, choosing $\theta$ close enough to $\frac{d\theta_1}{2}$, we see that
$(1-\theta)(1+(\alpha-1)\theta_1)>1$ and upon taking $\delta$ small enough we obtain
$$
\|(u-\delta)_-(t)\|_{\infty} \leq \frac{\delta}{2},
$$
for $t\geq t_0>0$. The Lemma is proved. \qed
\section{Convergence to the semi-trivial solution $(0,\theta_\mu)$}
Through this Section we additionally assume that there exist constants $0<c_m<C_M$ and $\alpha \geq 1$ such that
\begin{equation}\label{hthetamu}
c_m s^\alpha \leq V(s) \leq C_M s^\alpha \; \mbox{ for all $s\in [0,\|u\|_{\infty}]$}.
\end{equation}
\begin{remark}
Let us observe that when $V'(0) \not = 0$ and (\ref{hipotesisV}) holds, then (\ref{hthetamu}) is true for $\alpha=1$. Moreover if $V\in \go C^k$ for $k\geq 1$ with $V^k (0) \not =0$ and $V^j (0)=0$ for $j<k$, then (\ref{hthetamu}) holds
true for $\alpha=k$. 
\end{remark}
In the following Theorem, we eliminate the restriction on $\mu$ of Theorem \ref{convul0}. However, we require the additional condition (\ref{hthetamu}) on $V$.
\begin{theorem}\label{thm4.2}
Let $\lambda =0$ and assume (\ref{hthetamu}), then
$$
\lim_{t \rightarrow +\infty} \| u(t)\|_{W^{m,p}} =0,
$$
for any $m<1$ and $p\geq 2$.
\end{theorem}
\noindent {\bf Proof.} On the one hand, we multiply the $u$-equation of (\ref{sis00})
by $u$ and we integrate in the space variable to obtain
\begin{equation} \label{sss1.2}
\begin{array}{ll}
\displaystyle{ \frac{d}{2dt} \int_\Omega u^2 } &=\displaystyle{ \int_\Omega \left ( -|\nabla u|^2 + V(u)\nabla v \cdot \nabla u -u^3\right )- \mu\int_{\Gamma_2} \frac{V(u)uv}{1+v}} \\
&= \displaystyle{\int_\Omega \left ( -|\nabla u|^2 + \nabla v \cdot \nabla \varphi (u) -u^3\right )- \mu\int_{\Gamma_2} \frac{V(u)uv}{1+v}},
\end{array}
\end{equation}
with
$$
\varphi (u)=\int_0^u V(s)\,ds.
$$
On the other hand, we multiply the $v$-equation of (\ref{sis00}) by $\varphi (u)$. Integrating in space, we obtain
$$
\int_\Omega \nabla v \cdot \nabla \varphi (u) = -\int_\Omega \varphi (u) v_t + \mu \int_{\Gamma_2} \frac{v\varphi (u)}{1+v} -\int_\Omega v \varphi (u) - \int_\Omega c uv\varphi (u).
$$
Inserting the above equality into (\ref{sss1.2}) we have
\begin{equation} \label{sss2}
\displaystyle{ \frac{d}{2dt} \int_\Omega u^2 } =\displaystyle{ \int_\Omega \left ( -|\nabla u|^2 + \varphi (u) (-v_t -v -cuv)  -u^3\right )+ \mu\int_{\Gamma_2} \frac{v}{1+v}(\varphi (u)-V(u)u)}.
\end{equation}
Next we estimate $v_t$. Multiplying the $v$-equation by $v_t$ and integrating over $\Omega$ we see that
$$
\frac{1}{2} \int_\Omega v_t^2 + \frac{d}{2 dt}\int_\Omega |\nabla v|^2 + \frac{d}{2dt}\int_\Omega v^2 - \frac{\mu d}{dt} \int_{\Gamma_2} \theta (v) = -\int_\Omega c uvv_t\,,
$$
where
$$
\theta (v) := \int_0^v \frac{s}{1+s} ds.
$$
Therefore, by the uniform bound of $v$ in $\mathcal C(\overline \Omega)$ we deduce
$$
\frac{1}{4} \int_\Omega v_t^2 +  \frac{d}{2 dt}\int_\Omega |\nabla v|^2 + \frac{d}{2dt}\int_\Omega v^2 - \frac{\mu d}{dt} \int_{\Gamma_2} \theta (v) \leq M \int_\Omega u^2.
$$
After integrating over the interval $(\tau,t)$ we find, by Lemma \ref{lforu}, that for $t \geq \tau$
\begin{equation} \label{nwe}
\int_\tau^t \int_\Omega v_t^2 \leq C.
\end{equation}
Next, by (\ref{hthetamu}) we obtain from (\ref{sss2}) that
\begin{equation} \label{g1h1}
\begin{array}{ll}
\displaystyle{ \frac{d}{2dt} \int_\Omega u^2} &\leq \displaystyle{-\int_\Omega |\nabla u|^2 + \int_\Omega \varphi (u)^2 + \int_\Omega v_t^2 + \mu \int_{\Gamma_2} \frac{v C_u u^{\alpha+1}}{(1+v)(\alpha+1)}} \\
& \leq \displaystyle{-\int_\Omega |\nabla u|^2 + \max_{s \in [0,C_u]} V^2 (s) \int_\Omega u^2 + \int_\Omega v_t^2 + \frac{\mu C_u^2}{\alpha +1} \int_{\Gamma_2} \frac{v u^\alpha}{1+v}}.
\end{array}
\end{equation}
By Lemma \ref{lforu} and (\ref{hthetamu}) we get
\begin{equation}\label{raz}
\int_\tau^t \int_{\Gamma_2} \frac{v u^\alpha}{1+v} \leq C
\end{equation}
for $t \geq \tau$.
According to (\ref{raz}) and (\ref{nwe}) we find upon integration of (\ref{g1h1}) over the time interval $(\tau, t)$ that for $t \geq \tau$
$$
\int_\tau^t \int_\Omega |\nabla u|^2 \leq C.
$$
From the last estimate, a similar argument to the one used previously yields
$$
\int_\tau^t \left | \frac{d}{dt} \int_\Omega u^2\right | \leq C
$$
for $t \geq \tau$.
Thus, by Lemma \ref{convLema}
$$
\lim_{t \rightarrow +\infty} \| u(t) \|_2 =0.
$$
Finally, we can infer the result arguing as in the end of the proof of Theorem~\ref{convul0}. \qed

Next we prove a lemma which we will use in the proof of Theorem \ref{thm4.4}. As a by-product of the following 
lemma we learn a qualitative information that $v$ is bounded away from $0$ for times large enough. We shall obtain a lower bound on $v$ by considering a subsolution to an elliptic problem which is also a subsolution to a second equation in (\ref{sis00}). 
\begin{lemma}\label{lem4.3}
Let $\lambda =0$ and $\mu > \mu_1$. If the condition (\ref{hthetamu}) is satisfied then
there exist constants $c_1, \tau >0$ such that for $t \geq \tau$
\begin{equation} \label{estvl}
v(t) >c_1.
\end{equation}
\end{lemma}
\noindent {\bf Proof.} Let $k \in (\mu_1,\mu)$. Since $\lambda_1 (-\Delta + 1;{\go N};{\go N} + b(x))$ is increasing with respect to $b$ (see \cite[Proposition 3.3]{cclp}), we have  
$$
\lambda_1 (-\Delta + 1; {\go N};{\go N} -\mu) <\lambda_1 (-\Delta + 1; {\go N}; {\go N} -k) < \lambda_1 (-\Delta + 1; {\go N};{\go N} -\mu_1 )=0.
$$
Therefore, there exists $\epsilon >0$ such that
$$
\lambda_1 (-\Delta +1 ; {\go N};{\go N} - k)= -c\epsilon_0 \mbox{ i.e. $\lambda_1 (-\Delta +1 +c\epsilon_0; {\go N};{\go N}-k)=0$. }
$$
Let $\varphi_1$ be the positive eigenfunction with $\|\varphi_1\|_\infty =1$ associated to the above eigenvalue i.e. $\varphi_1$ satisfies 
$$
\left \{
\begin{array}{ll}
-\Delta \varphi_1 + (1+\epsilon_0 c)\varphi_1 = 0 & \mbox{in $\Omega$,}\\
\displaystyle\frac{\p \varphi_1}{\p n} =0 & \mbox{on $\Gamma_1$,}\\
\displaystyle\frac{\p \varphi_1}{\p n} = k \varphi_1 & \mbox{on $\Gamma_2$.}
\end{array}
\right.
$$
By Theorem \ref{thm4.2} there exists $t_0 >0$ such that  $0\leq u(t) < \epsilon_0$ for all $t\geq t_0>0$. We claim that there exists $\delta >0$ such that $\underline w = \delta \varphi_1 $ is a subsolution to
$$
\left \{
\begin{array}{ll}
w_t-\Delta w + (1+cu)w = 0 & \mbox{in $\Omega \times (t_0,+\infty)$,}\\
\displaystyle\frac{\p w}{\p n} =0 & \mbox{on $\Gamma_1\times (t_0,+\infty)$,}\\
\displaystyle\frac{\p w}{\p n} = \mu \frac{w}{1+w} & \mbox{on $\Gamma_2\times (t_0,+\infty)$.}\\
w(x,t_0)=v(x,t_0) & \mbox{in $\Omega$.}\\
\end{array}
\right.
$$
Therefore $v(x,t) \geq \delta \varphi_1 \geq c_1$. It remains to prove the claim. By the strong maximum principle $v(x,t_0) > c>0$. Thus there exists $\delta >0$ such that $\delta \varphi_1 < v(x,t_0)$. Moreover, choosing $\delta >0$ such that $k(1+\delta)<\mu$ we make sure that  
\[
\frac{\p \underline w}{\p n} \leq \mu \frac{\underline w}{1+ \underline w} 
\]
on $\Gamma_2\times (t_0,+\infty)$. Hence the claim is shown and the lemma follows.
\qed
Now we are in a position to prove the main result of this section. To this end we make use of the theorem by Amann and L\'opez-G\'omez, see \cite{AJ}, stating the equivalence between positivity of principal eigenvalue and existence of stricly positive supersolution of some elliptic problems (the previous version of this theorem for the Dirichlet problem was shown in \cite{LG}). 

\begin{theorem}\label{thm4.4}
Let $\lambda =0$ and assume (\ref{hthetamu}), then
$$
\lim_{t \rightarrow +\infty} \| v(t)- \theta_\mu \|_{2} =0.
$$
\end{theorem}
\noindent {\bf Proof.} Let $z(t)= v(t) -\theta_\mu$. Then $z$ solves the following parabolic problem
\begin{equation} \label{sisz} \left\{
\begin{array}{ll}
z_t =\Delta z -z -cuv &\mbox{in $\O\times (0,T)$,}\\
\displaystyle\frac{\p v}{\p n}= 0 &\mbox{on $\G_1\times (0,T)$,}\\
\displaystyle\frac{\p z}{\p n}=
\m\displaystyle\frac{z}{(1+v)(1+\theta_\mu)} &\mbox{on $\G_2\times (0,T)$,}\\
z(x,0)=v_0(x)-\theta_\mu  &\mbox{in $\O$.}
\end{array}\right.
\end{equation}
We multiply (\ref{sisz}) by $z$ to obtain
\begin{equation}\label{zini}
\frac{d}{2dt} \int_\Omega z^2 = -\int_\Omega |\nabla z|^2 + \mu \int_{\Gamma_2} \frac{z^2}{(1+v)(1+\theta_\mu)} - \int_\Omega z^2 - \int_\Omega c uvz.
\end{equation}
In order to estimate the right-hand side of (\ref{zini}) for $t \geq t_0$, we pick $\gamma >1$ such that 
\begin{equation}\label{c_1}
\frac{\gamma}{1+ c_1}<1
\end{equation} 
where $c_1$ is given in (\ref{estvl}). For each $t \geq t_0$ we consider the eigenvalue problem
\begin{equation} \label{sisw}
\left\{
\begin{array}{ll}
-\Delta w +w = \lambda w &\mbox{in $\O$,}\\
\displaystyle\frac{\p w}{\p n}= 0 &\mbox{on $\G_1$,}\\
\displaystyle\frac{\p w}{\p n}=
\displaystyle\frac{\mu \gamma w}{(1+v(t))(1+\theta_\mu)} &\mbox{on $\G_2$.}
\end{array}\right.
\end{equation}
Next, we see that $\theta_\mu$ is a strict supersolution of
$$
\left\{
\begin{array}{ll}
-\Delta w +w = 0 &\mbox{in $\O$,}\\
\displaystyle\frac{\p w}{\p n}= 0 &\mbox{on $\G_1$,}\\
\displaystyle\frac{\p w}{\p n}=
\displaystyle\frac{\mu \gamma w}{(1+v(t))(1+\theta_\mu)} &\mbox{on $\G_2$.}
\end{array}\right.
$$
Indeed,
$$
-\Delta \theta_\mu + \theta_\mu =0 \mbox{ in $\Omega$,}
$$
$$
\frac{\partial \theta_\mu}{\partial n} =0 \mbox{ on $\Gamma_1$,}
$$
Finally by the choice of $\gamma$ (see (\ref{c_1})) and Lemma \ref{lem4.3} we have
$$
\frac{\partial \theta_\mu}{\partial n} = \mu\frac{ \theta_\mu}{1 + \theta_\mu} > \frac{\mu \gamma \theta_\mu}{(1+v(t))(1+\theta_\mu)} \mbox{ on $\Gamma_2$.}
$$
Therefore, by \cite[Theorem 2.4]{AJ} we get $\lambda_1 >0$, the principal eigenvalue of (\ref{sisw}). Next, the variational characterization of the principal eigenvalue entails
$$
\lambda_1 = \inf_{\varphi \in H^1 (\Omega)} \frac{ \displaystyle{\int_\Omega |\nabla \varphi |^2 + \int_\Omega \varphi^2 - \mu \gamma \int_{\Gamma_2} \frac{\varphi^2}{(1+v(t))(1+\theta_\mu)}}}{\displaystyle\int_\Omega \varphi^2 }\,.
$$
Thus, for all $\varphi \in H^1( \Omega)$ we have
$$
\lambda_1 \gamma^{-1} \int_\Omega \varphi^2 \leq \gamma^{-1} \int_\Omega |\nabla \varphi|^2 + \gamma^{-1} \int_\Omega \varphi^2 - \mu \int_\Omega \frac{\varphi^2}{(1+v(t))(1+\theta_\mu)} \,.
$$
In particular, we can apply it in (\ref{zini}) to obtain the following inequality
$$
\frac{d}{2dt} \int_\Omega z^2 + (1-\gamma^{-1}) \int_\Omega |\nabla z|^2 + (1-\gamma^{-1} + \lambda_1 \gamma^{-1}) \int_\Omega z^2 \leq \int_\Omega cuvz \,.
$$
Therefore, there exists $M>0$ such that
$$
\frac{d}{2dt} \int_\Omega z^2 + (1-\gamma^{-1}) \left ( \int_\Omega |\nabla z|^2 + \int_\Omega z^2 \right ) \leq M \int_\Omega u^2 \,.
$$
Integrating the above estimate on the time interval $(\tau,t)$ we obtain for $t \geq \tau$,
\begin{equation}\label{dwa}
\int_\tau^t \int_\Omega |\nabla z|^2 + \int_\Omega z^2 \leq C.
\end{equation}
In view of (\ref{dwa}) one infers
$$
\int_\tau^t \left | \frac{d}{dt} \int_\Omega z^2\right | \leq C
$$
for $t \geq \tau$. Finally, the result follows  by Lemma \ref{convLema}. \qed

\noindent
{\bf Acknowledgement.}T. Cie\'slak was partially supported by the Polish Ministry of Science and Higher Education under grant number NN 201366937. C. Morales-Rodrigo was supported by Ministerio de Ciencia e Innovaci\'on and FEDER under grant MTM2009-12367.

\end{document}